\documentclass[11pt]{article}

\usepackage[a4paper]{geometry}
\usepackage{authblk}
\usepackage[utf8]{inputenc} 
\usepackage[T1]{fontenc}    
\usepackage{hyperref}       

\usepackage[fleqn]{amsmath}
\usepackage{amsfonts}
\usepackage{amsthm}
\usepackage{amsbsy}
\usepackage{bbm}
\usepackage{enumitem}
\usepackage{graphicx}
\usepackage{mathtools,siunitx}
\usepackage{soul} 

\newtheorem{defn}{Definition}
\newtheorem{prop}{Proposition}
\newtheorem{cor}{Corollary}
\newtheorem{lem}{Lemma}

\newtheorem{thm}{Theorem}

\theoremstyle{remark}
\newtheorem*{rem}{Remark}

\newcommand{\me}{\mathrm{e}}
\newcommand{\Prob}{\mathbb{P}}
\renewcommand{\d}{\mathrm{d}}
\newcommand{\mat}[1]{\mathbf{#1}}
\newcommand{\Rset}{\mathbb{R}}

\begin{document}

\title{The Markov chain embedding problem in a low jump frequency context}

\author{Philippe Carette}
\affil{Department of Economics\\ Ghent University, Sint-Pietersplein 5, B-9000 Ghent, Belgium \\
\texttt{philippe.carette@ugent.be}}

\author{Marie-Anne Guerry}
\affil{Department Business Technology and Operations\\ Vrije Universiteit Brussel, Pleinlaan 2, B-1050 Brussels, Belgium\\
\texttt{marie-anne.guerry@vub.be}}
 
\maketitle

\begin{abstract}
We consider the problem of finding the transition rates of a continuous-time homogeneous Markov chain under the empirical condition that the state changes at most once during a time interval of unit length. It is proven that this conditional embedding approach results in a unique intensity matrix for a transition matrix with non-zero diagonal entries. Hence, the presented conditional embedding approach has the merit to avoid the identification phase as well as regularization for the embedding problem. 
The resulting intensity matrix is compared to the approximation for the Markov generator found by Jarrow in \cite{jarrow1997markov}.

\end{abstract}

\section{Introduction}

The embedding problem of Markov chains is a long standing problem where a given stochastic matrix is examined as the 1-step transition matrix of some continuous-time homogeneous Markov chain (CTHMC) (\cite{elfving,kingman}). This problem boils down to characterizing the empirical transition matrix $\mat{\widehat{P}}$ as the exponential of some matrix $\mat{Q}$ with all non-negative off-diagonal entries and zero row-sums, called an intensity matrix. This matrix $\mat{Q}$ represents the transition rates of the underlying CTHMC. If such a $\mat{Q}$ exists, $\mat{\widehat{P}}$ is said to be embeddable. It turns out that the embedding problem is a formidable one in a number of respects. First, $\mat{\widehat{P}}$ may not be embeddable. In that case, a regularization algorithm can be used to find an intensity matrix $\mat{Q}$ for which $||\mat{\widehat{P}} - \exp(\mat{Q})||$ is minimized (\cite{davies2010embeddable,kreinin2001regularization, israel2001finding}). Next, no embeddability criteria in terms of the matrix elements, which are easily verifiable in practice, seem at hand when the number of states exceeds 3. Lastly, for an embeddable $\mat{\widehat{P}}$, there may not be a unique solution to the equation $\exp(\mat{Q}) = \mat{\widehat{P}}$ in the set of intensity matrices. The identification aspect of the embedding problem deals with the selection of the suitable intensity matrix reflecting the nature of the system under study (\cite{singer1974spilerman}).

More recently, model specific embedding problems are studied for specific subcategories of transition matrices (\cite{ardiyansyah2021model,baake2020notes, roca2018embeddability, jarrow1997markov}). In modeling a specific context, the transition matrix as well as the generator matrix are expected to reflect the characteristics of the system under study. Hence, the transition matrix is subjected to constraints and, therefore, belongs to a specific subset of stochastic matrices, and similar, the generator matrix is expected to be an element of a specific subset of intensity matrices. Whereas model specific embedding problems are characterized by setting model assumptions and restrictions on the transition matrix, this paper presents an embedding approach that incorporates empirical assumptions.

More specifically, we propose the conditional embedding approach where the empirical 1-step transition matrix $\mat{\widehat{P}}$ corresponds with the conditional 1-step transition matrix of the CTHMC given the event that at most one jump has occurred during a time interval of unit length. For a Markov model the unit time interval can be defined in such a way that the empirical 1-step transition matrix meets this condition. Moreover, this condition is inherent in some applications. For example, in credit rating migration models the credit ratings are typically based on slowly varying characteristics, such that they do not tend to change more than once within the baseline time interval (e.g. a quarter).

We found that, regardless the number of states, exactly one intensity matrix solves this conditional embedding problem when $p_{ii} > 0$ for all $i$. Our approach results in an easy embeddability criterium and does not require identification neither regularization. Moreover, the presented conditional embedding approach and its proven properties, result in an embeddability roadmap reflecting that the conditional embedding approach is atmost useful in case either the transition matrix is not embeddable or no unique Markov generator can be identified based on the context of the system.

\bigskip




\section{Conditional transition probabilities}

In order to state the conditional embedding problem, we first introduce the concept of conditional transition probability.

Consider a continuous-time homogeneous Markov chain (CTHMC) $(X_t)_{t\geq 0}$ on a probability triple $(\Omega,\mathcal{F},\Prob)$ with state space $\mathcal{S}=\{1,2,\dots,n\}$.

\begin{defn}
Let $E\in\mathcal{F}$. We call the matrix $\mat{P}^{E}$ with elements 
\[
p^{E}_{ij}=\Prob(X_1=j\,|\,X_0=i\,,\,E),\quad i,j\in\mathcal{S}\text{,}
\]
the \emph{conditional one-step transition probability matrix given the event $E$} of the chain $(X_t)_{t\geq0}$.
\end{defn}

\begin{rem}
    The usual (unconditional) one-step transition probabilities $p_{ij}=\Prob(X_1=j\,|\,X_0=i)$ can be obtained by setting $E = \Omega$, that is, $p_{ij}=p^{\Omega}_{ij}$.
\end{rem}

In the remainder of this paper, we are interested in the event $E=\{N_{\!J}\leq 1\}$, where $N_{\!J}$ is the random variable counting the  state changes or jumps of the CTHMC up to time $t=1$.

The relationship between the conditional transition matrix $\mat{P}^{\{N_{\!J}\leq 1\}}$ and the transition rate matrix $\mat{Q}$ of the CTHMC is given by the following proposition.

\begin{prop} \label{prop:ptilde_Q}
 For a CTHMC with transition rate matrix $\mat{Q}=(q_{ij})$, it holds that 
 \[
 {p}^{\{N_{\!J}\leq 1\}}_{ij} = \frac{p_{ij}^{*}}{\sum_{k=1}^{n} p_{ik}^{*}}\quad\text{for all $i$ and $j$,}
 \]
 where
 \[
 p_{ij}^{*} = \begin{cases}
q_{ij}\,\tau(q_{ii},q_{jj}) & \text{if $i\neq j$}\\
\tau(q_{ii},q_{ii}) & \text{if $i=j$}
\end{cases}
 \]
and where the function $\tau:\Rset^2\to\Rset$ is defined as
\begin{equation}\label{eq:tau}
\tau(x,y)=\int_{0}^{1} \me^{ux+(1-u)y}\,\d{u}
=
\begin{cases}
\frac{\me^{x}-\me^{y}}{x-y} & \text{if $x\neq y$}\\
\me^{x} & \text{if $x=y$}
\end{cases}
\text{.}
\end{equation}
\end{prop}

\begin{proof}
Using the definition of conditional probability, 
\[
\Prob(A\,|\,B\cap C)=\frac{\Prob(A\cap B\,|\,C)}{\Prob(B\,|\,C)},\quad\text{if $\Prob(B\,|\, C)>0$.}
\]
 Hence, 
\begin{equation}\label{eq:ptildeij}
{p}^{\{N_t\leq 1\}}_{ij} = \Prob(X_1=j\,|\,X_0=i,N_{\!J}\leq 1)=\frac{\Prob(X_1=j,N_{\!J}\leq 1\,|\,X_0=i)}{\Prob(N_{\!J}\leq 1\,|\,X_0=i)}\text{.}
\end{equation}
Let us denote $p^{*}_{ij}=\Prob(X_1=j,N_{\!J}\leq 1\,|\,X_0=i)$. Using the sum rule for disjoint events, we then have 
\[
{p}^{\{N_{\!J}\leq 1\}}_{ij}=\frac{p_{ij}^{*}}{\sum_{k=1}^{n} p_{ik}^{*}}\text{.}
\]

Let us now calculate $p^\ast_{ij}$, which is the joint probability of being in state $j$ at $t=1$ in at most one jump, starting from state $i$ at $t=0$.  For $k\in\{1,\dots,n\}$, let $f_{k}$ be the density function of the holding time $H_k$ in state $k$ and $F_k$ the associated cumulative distribution function. 
For $i\neq j$, denote by $s_{ij}$ the transition probability from state $i$ to state $j$ conditional on transitioning out of state $i$. Marginalising on $H_i$, we find for $i\neq j$
\begin{align*}
p^{*}_{ij}
&=\Prob(H_i<1,X_{H_i}=j,H_j>1-H_i\,|\,X_0=i)\\
&=\int_{0}^{1}\Prob(X_u=j,H_j>1-u\,|\,H_i=u,X_0=i)f_i(u)\,\d{u} \\
&=\int_{0}^{1}\Prob(H_j>1-u\,|\,X_u=j,H_i=u,X_0=i)\Prob(X_u=j\,|\,H_i=u,X_0=i)f_i(u)\,\d{u} \\
&=\int_{0}^{1}(1-F_j(1-u))\,s_{ij}\,f_i(u)\,\d{u}\qquad 
\end{align*}
and 
\[
p^{*}_{ii}=\Prob(H_i>1\,|\,X_0=i)=1-F_i(1)\text{.}
\]

Since $H_k$ has an exponential distribution with rate parameter $-q_{kk}$ and since $s_{ij}=\frac{q_{ij}}{-q_{ii}}$ ($i\neq j$), we have for $i\neq j$
\[
p^{*}_{ij}=\int_{0}^{1} \me^{q_{jj}(1-u)} \,\frac{q_{ij}}{-q_{ii}}\,(-q_{ii})\me^{q_{ii}u}\,\d{u}
= q_{ij}\int_{0}^{1}\me^{q_{ii}u+q_{jj}(1-u)}\,\d{u}
=q_{ij}\,\tau(q_{ii},q_{jj})
\]
where the function $\tau$ is defined as in (\ref{eq:tau}). Finally, $p^{*}_{ii}=1-F_i(1)=\me^{q_{ii}}=\tau(q_{ii},q_{ii})$.
\end{proof}

One can remark that Minin et al. \cite{minin2008counting} arrive at the same result for $p^\ast_{ij}$ using a recursive relation for the joint probabilities $\Prob(X_1=j\,|\,X_0=i,N_{\!J}= n)$, $i,j\in\mathcal{S}$.  

\begin{cor} \label{diagonal elements strict pos} For all $i$, we have that ${p}^{\{N_{\!J}\leq 1\}}_{ii}>0$.
\end{cor}
\begin{proof}
    This follows from the fact that $p^{*}_{ii}=\tau(q_{ii},q_{ii})=\me^{q_{ii}}>0$. 
\end{proof}

\begin{rem}
An alternative argument for corollary~1 goes as follows. Given the event $\{N_{\!J}\leq 1\}$, the only way of going from state $i$ at time $t=0$ to state $i$ at time $t=1$, is to remain in that state throughout the entire time interval from $t=0$ to $t=1$. The probability of this event is non-zero, since the holding time in a state has an exponential distribution.
\end{rem}

    According to proposition~\ref{prop:ptilde_Q}, the conditional transition matrix $\mat{P}^{\{N_{\!J}\leq 1\}}$ depends on the transition rate matrix $\mat{Q}$ of the CTHMC involved. In what follows, and when needed, we explicitly indicate this dependency using the notation $\mat{P}^{\{N_{\!J}\leq 1\}}(\mat{Q})$.



\section{Conditional embedding problem}

When building a discrete time Markov model, the choice of the time unit and time interval is important to end up with a valid model (\cite{Bhuller}). In this respect, an appropriate choice can be made by comparing for diverse values of the time unit the internal validity of the corresponding models. The internal validity of a model is determined by the discrepancy between the observed stock vectors and the stock vectors that are estimated by the model. Based on goodness of fit tests a time unit can be selected that results in a model for which the discrepancy between observed and estimated stock vectors is limited.  (\cite{Sales}). For an appropriate time unit it is acceptable to assume that there is at most $1$ jump in between $t=0$ and $t=1$. Indeed, more than $1$ jump during a one-unit time interval would result in a situation where the transitions to and from the intermediate state are not captured by the discrete time Markov model.

A question that then naturally arises is whether, for a given stochastic matrix $\mat{P}$, there does exist an intensity matrix $\mat{Q}$ such that $\mat{P}^{\{N_{\!J}\leq 1\}}(\mat{Q})=\mat{P}$. And if so, whether such an intensity matrix $\mat{Q}$ is unique.

It will be helpful to introduce some terminology before proceeding.

\begin{defn}
    A stochastic matrix $\mat{P}$ is called \emph{$J_1$-embeddable} iff there exist a CTHMC with transition rate matrix $\mat{Q}$ satisfying $\mat{P}=\mat{P}^{E}(\mat{Q})$, where $E$ is the event that the CTHMC changes state at most once between $t=0$ and $t=1$. Such a transition rate matrix is called a \emph{$J_1$-generator} of $\mat{P}$.
\end{defn}

For a transition matrix $\mat{P}$ that is not embeddable a \emph{$J_1$-generator} can be seen as a solution to the generalization problem where the intensity matrix $\mat{Q}$ satisfies  $\mat{P}^{\{N_{\!J}\leq 1\}}(\mat{Q})=\mat{P}$. For a transition matrix $\mat{P}$ that is embeddable, its Markov generator not necessarily equals the \emph{$J_1$-generator}. In fact, the solution to the conditional embedding problem is generally different from the solution to the (general) embedding problem if the latter exists. However, if the time-scale is chosen such that no more than one transition occurs in the system during the unit time interval, we might expect Markov generator to be close to \emph{$J_1$-generator} in some sense.

A matrix $\mat{P}$ that is embeddable satisfies the necessary condition for embeddability formulated by Goodman in \cite{Goodman}: $\prod_{i=1}^{n} p_{ii} \geq \det{P} > 0$.
Hence, all diagonal entries of such an embeddable matrix $\mat{P}$ are non-zero. Consequently, a matrix $\mat{P}$ with $p_{ii} = 0$, for some $i$, is neither embeddable nor $J_1$-embeddable, according to corollary \ref{diagonal elements strict pos}. For this reason, and without loss of generality, we examine in the remainder of the paper stochastic matrices $\mat{P}=(p_{ij})$ satisfying $p_{ii}>0$ for all $i$. 

It turns out that the off-diagonal elements of a $J_1$-generator of $\mat{P}$ are uniquely determined by its diagonal elements and the elements of $\mat{P}$. To formulate this relationship, we introduce the  function $\rho:\Rset_+^2\to\Rset_+$, defined as follows:
\begin{equation}\label{eq:rho(x,y)}
    \rho(x,y)=\frac{\me}{\tau(1-\ln{x},1-\ln{y})}
    =\begin{cases}
         xy\frac{\ln{x}-\ln{y}}{x-y} & \text{if $x\neq y$}\\
         x & \text{if $x=y$.}
        \end{cases} 
\end{equation}

\begin{prop}\label{prop:q_ij}
    Suppose $\mat{P}=(p_{ij})$ is a $n\times n$ stochastic matrix satisfying $p_{ii}>0$ for all $i$. If $\mat{Q}=(q_{ij})$ is a $J_1$-generator of $\mat{P}$, then
    \begin{equation*}
	   q_{ij}=
    \frac{\rho(\theta_i,\theta_j)p_{ij}}{\theta_ip_{ii}}
    ,\quad \text{for all $i\neq j$,}
    \end{equation*}
    where $\theta_i=\me^{1-q_{ii}}$ for all $i$ and the function $\rho:\Rset_+^2\to\Rset_+$ is given by \eqref{eq:rho(x,y)}.
\end{prop}
\begin{proof}
    Suppose $\mat{Q}$ is a $J_1$-generator of $\mat{P}=(p_{ij})$ and let $i\neq j$. Then, according to proposition~\ref{prop:ptilde_Q}, we have
    \[
    \frac{p_{ij}}{p_{ii}}=\frac{{p}^{\{N_{\!J}\leq 1\}}_{ij}}{{p}^{\{N_{\!J}\leq 1\}}_{ii}}=\frac{p^*_{ij}}{p^*_{ii}}=\frac{q_{ij}\tau(q_{ii},q_{jj})}{\tau(q_{ii},q_{ii})}\text{.}
    \]
    Consequently, since $\tau(q_{ii},q_{ii})=\me^{q_{ii}}=\me/\theta_i$ and $\tau(q_{ii},q_{jj})=\tau(1-\ln{\theta_i},1-\ln{\theta_j})=\me/\rho(\theta_i,\theta_j)$, we get
    \[
    q_{ij}=\frac{\tau(q_{ii},q_{ii})p_{ij}}{\tau(q_{ii},q_{jj})p_{ii}}=\frac{(\me/\theta_i)p_{ij}}{(\me/\rho(\theta_i,\theta_j))p_{ii}}=\frac{\rho(\theta_i,\theta_j)p_{ij}}{\theta_i p_{ii}},\quad \text{for all $i\neq j$.}\qedhere
    \]
\end{proof}

The result of proposition~\ref{prop:cembed:Tfixpt} yields a condition on the diagonal elements of any $J_1$-generator of $\mat{P}$.

\begin{prop}\label{prop:cembed:Tfixpt}
    Suppose $\mat{P}=(p_{ij})$ is a $n\times n$ stochastic matrix satisfying $p_{ii}>0$ for all $i$. Then, if $\mat{Q}=(q_{ij})$ is a $J_1$-generator of $\mat{P}$, the $n$-tuple $(\me^{1-q_{11}},\dots,\me^{1-q_{nn}})$ is a fixed point of the vector function $\mat{T}=(T_1,\dots,T_n):\Rset_+^{n}\to\Rset_+^{n}$ defined as follows
    \begin{equation}\label{eq:T}
        T_i(x_1,\dots,x_n)=\exp{W_0\Bigl(\frac{1}{p_{ii}}\sum_{j=1}^{n}p_{ij}\rho(x_i,x_j)\Bigr)}\quad\text{for all $i$,}
    \end{equation}
    where $W_0$ denotes the principal branch of the Lambert W function and where the function $\rho:\Rset_+^2\to\Rset_+$ is defined as in \eqref{eq:rho(x,y)}.
\end{prop}
\begin{proof}
    Denote $\theta_i=\me^{1-q_{ii}}$ for all $i$. Then $\theta_i>0$ and $q_{ii}=1-\ln{\theta_i}$ for all $i$.  Using proposition~\ref{prop:q_ij} and the fact that $\mat{Q}$ is an intensity matrix, we then have
    \[
    -1+\ln{\theta_i}=-q_{ii}=\sum_{j:j\neq i}q_{ij}=\sum_{j:j\neq i}\frac{\rho(\theta_i,\theta_j)p_{ij}}{\theta_i p_{ii}},\quad\text{for all $i$,}
    \]
    which can be rewritten, using the fact that $\rho(\theta_i,\theta_i)=\theta_i$, as
    \begin{equation}\label{eq:theta*ln(theta)}
        \theta_i\ln{\theta_i}=\frac{1}{p_{ii}}\sum_{j=1}^{n}p_{ij}\rho(\theta_i,\theta_j)\quad\text{for all $i$.}
    \end{equation}
    Using the principal branch $W_0$ of the Lambert W function (which is the multi-valued inverse of the function $w\mapsto w\me^{w}$ ($w\in\mathbbmss{C}$), see \cite{corless1996lambertw}), we find that 
    \[
    \ln{\theta_i}=W_0\Bigl(\frac{1}{p_{ii}}\sum_{j=1}^{n}p_{ij}\rho(\theta_i,\theta_j)\Bigr)\quad\text{for all $i$,}
    \]
    which proves the result.
\end{proof}

Proposition~\ref{prop:q_ij} and proposition~\ref{prop:cembed:Tfixpt} entail that a $J_1$-generator of $\mat{P}$ defines a fixed point of the vector function $\mat{T}$. The converse is also true, as is stated in the following proposition.
\begin{prop}\label{prop:fixpT_to_J1gen}
    Let the stochastic $n\times n$ matrix $\mat{P}=(p_{ij})$ be such that $p_{ii}>0$ for all $i$. Suppose $\mat{\theta}=(\theta_1,\dots,\theta_n)$ is a fixed point of the vector function $\mat{T}:\Rset_+^n\to\Rset_+^n$, defined  in \eqref{eq:T}. Then, the matrix $\mat{Q}=(q_{ij})$ with elements 
    \begin{equation*}
         q_{ii}= 1-\ln{\theta_i},\qquad q_{ij}=\frac{\rho(\theta_i,\theta_j)p_{ij}}{\theta_i p_{ii}}\quad(i\neq j)
    \end{equation*}
    where $\rho$ is defined by \eqref{eq:rho(x,y)}, is a $J_1$-generator of $\mat{P}$. 
\end{prop}
\begin{proof}
    Let $\mat{\theta}=(\theta_1,\theta_2,\dots,\theta_n)\in\Rset_+^n$ be a fixed point of $\mat{T}$ and let the matrix $\mat{Q}$ be constructed as stated above. We first show that $\mat{Q}$ is an intensity matrix.  By definition, all off-diagonal elements of $\mat{Q}$ are non-negative. Since $T_i(\mat{\theta})=\theta_i$, we have $\ln{\theta_i}=W_0(\frac{1}{p_{ii}}\sum_{j}p_{ij}\rho(\theta_i,\theta_j))$ yielding $\theta_i\ln{\theta_i}=\frac{1}{p_{ii}}\sum_{j}p_{ij}\rho(\theta_i,\theta_j)$, by definition of the Lambert $W_0$-function. Using $q_{ij}=\frac{\rho(\theta_i,\theta_j)p_{ij}}{\theta_i p_{ii}}$ ($i\neq j$) and $\rho(\theta_i,\theta_i)=\theta_i$, we can rewrite this equation as $\theta_i(1-q_{ii})=\theta_i+\frac{1}{p_{ii}}\sum_{j:j\neq i}q_{ij}\theta_ip_{ii}$. After simplification, we get $q_{ii}=-\sum_{j:j\neq i}q_{ij}$. Thus $\mat{Q}$ has zero row-sums. Consequently, $\mat{Q}$ is an intensity matrix. 
    
    It remains to be shown that ${p}^{\{N_{\!J}\leq 1\}}_{ij}(\mat{Q})=p_{ij}$ for all $i$ and $j$. By proposition~\ref{prop:ptilde_Q}, we have 
    ${p}^{\{N_{\!J}\leq 1\}}_{ij}(\mat{Q})=\frac{p^*_{ij}}{\sum_{k}p^*_{ik}}$, where $p^*_{ik}=q_{ik}\tau(q_{ii},q_{kk})$ if $i\neq k$ and $p^*_{ii}=\tau(q_{ii},q_{ii})$ and where the function $\tau:\Rset^2\to\Rset$ is defined by \eqref{eq:tau}. Using the definition of $\mat{Q}$ and \eqref{eq:rho(x,y)}, we then have that
    \[
    p^*_{ik}=\frac{\rho(\theta_i,\theta_k)p_{ik}}{\theta_i p_{ii}}\tau(1-\ln{\theta_i},1-\ln{\theta_k})=\frac{\me\, p_{ik}}{\theta_ip_{ii}},\qquad i\neq k
    \] 
    and 
    \[
    p^*_{ii}=\tau(q_{ii},q_{ii})=\me^{q_{ii}}=\me^{1-\ln{\theta_i}}=\frac{\me}{\theta_i}\text{.}
    \]
    Thus, $p^*_{ik}=\dfrac{\me\, p_{ik}}{\theta_ip_{ii}}$ for all $i$ and $k$, which yields
    \[
    \sum_{k}p^*_{ik}=\frac{\me}{\theta_i p_{ii}}\sum_k p_{ik}=\frac{\me}{\theta_i p_{ii}}\text{,}
    \]
    since $\mat{P}$ has all row sums equal to $1$. Consequently, for all $i$ and $j$, we get
    \[
    {p}^{\{N_{\!J}\leq 1\}}_{ij}(\mat{Q})=\frac{p^*_{ij}}{\sum_k p^*_{ik}}=\left.\frac{\me\,p_{ij}}{\theta_ip_{ii}}\middle/\frac{\me}{\theta_ip_{ii}}\right.=p_{ij},
    \]
    which concludes the proof.
\end{proof}

The following lemma states some properties of the vector function $\mat{T}$, which will play a crucial role in its number of fixed points.

\begin{lem}\label{lem:T}
    Let $\mat{P}=(p_{ij})$ be a $n\times n$ stochastic matrix. Let $\Delta=\max\{p_{11},\dots,p_{nn}\}$ and $\delta=\min\{p_{11},\dots,p_{nn}\}$. Suppose $\delta>0$. Consider the vector function $\mat{T}:\Rset_+^n\to\Rset_+^n$, defined as in \eqref{eq:T}, and the set
    \begin{equation}\label{eq:X}
        \mathcal{X}=\{\mat{x}=(x_1,\dots,x_n)\in\Rset^n\,|\, \forall i:\me^{1/\Delta}\leq x_i \leq \me^{1/\delta}\}\text{.}
    \end{equation}
    Then, 
    \begin{enumerate}[label=(\arabic*)]
        \item every fixed point of $\mat{T}$ belongs to $\mathcal{X}$. \label{lem:T:fp_in_X} 
        \item $\mat{T}$ maps $\mathcal{X}$ into $\mathcal{X}$. \label{lem:T:XtoX}
    \end{enumerate}
\end{lem}
\begin{proof} $ $
\begin{enumerate}[label=(\arabic*)] 
    \item 
    Let $\mat{\theta}=(\theta_1,\dots,\theta_n)\in\Rset_+^n$ be a fixed point of $\mat{T}$. Let $m=\min\{\theta_1,\dots,\theta_n\}$ and $M=\max\{\theta_1,\dots,\theta_n\}$. We shall prove that $m\geq\me^{1/\Delta}$ and that $M\leq\me^{1/\delta}$ . 
    
    Let $r$ be an index such that $\theta_r=m$. Then, by lemma~\ref{lem:rho}\ref{lem:rho:maxmin}, we have $\rho(\theta_r,\theta_j)\geq m$ for all $j$. Since $T_r(\mat{\theta})=\theta_r$, we have $\ln{\theta_r}=W_0(\frac{1}{p_{rr}}\sum_{j}p_{rj}\rho(\theta_r,\theta_j))$ yielding $\theta_r\ln{\theta_r}=\frac{1}{p_{rr}}\sum_{j}p_{rj}\rho(\theta_r,\theta_j)$ by definition of the Lambert $W_0$-function. Using the fact that $\sum_{j=1}^{n} p_{rj}=1$, we then obtain 
    \[
    m\ln{m}=\theta_r\ln{\theta_r}=\frac{1}{p_{rr}}\sum_{j=1}^{n}p_{rj}\rho(\theta_r,\theta_j)\geq \frac{m}{p_{rr}}\geq\frac{m}{\Delta}\text{,}
    \]
    which implies $\ln{m}\geq\frac{1}{\Delta}$ and thus $m\geq \me^{1/\Delta}$. To prove the second inequality, let $s$ be an index such that $\theta_s=M$. Then, by lemma~\ref{lem:rho}\ref{lem:rho:maxmin}, we have $\rho(\theta_s,\theta_j)\leq M$ for all $j$. Hence, by $T_s(\theta)=\theta_s$ and the unit row-sums property of $\mat{P}$, 
    \[
    M\ln{M}=\theta_s\ln{\theta_s}=\frac{1}{p_{ss}}\sum_{j=1}^{n}p_{sj}\rho(\theta_s,\theta_j)\leq \frac{M}{p_{ss}} \leq\frac{M}{\delta}\text{,}
    \]
    which yields $\ln{M}\leq \frac{1}{\delta}$ and thus $M\leq\me^{1/\delta}$.     
    \item 
    Let $(x_1,\dots,x_n)\in \mathcal{X}$. By lemma~\ref{lem:rho}\ref{lem:rho:maxmin},
    \begin{equation*}
        \me^{1/\Delta}\leq\min\{x_i,x_j\}\leq\rho(x_i,x_j)\leq\max\{x_i,x_j\}\leq \me^{1/\delta}\quad\text{for all $i$ and $j$.}
    \end{equation*}
    Then, since $\mat{P}$ has unit row sums, we have for all $i$
    \begin{equation*}
       \frac{\me^{1/\Delta}}{\Delta}\leq\frac{\me^{1/\Delta}}{p_{ii}}\leq \frac{1}{p_{ii}}\sum_{j=1}^np_{ij}\rho(x_i,x_j)\leq \frac{\me^{1/\delta}}{p_{ii}}\leq\frac{\me^{1/\delta}}{\delta}\text{.}    
    \end{equation*}
    Now, $W_0$ and $\exp$ are increasing functions, therefore
    \[
     \exp{W_0(\tfrac{1}{\Delta}\me^{1/\Delta})} \leq T_i(x_1,\dots,x_n)
    \leq\exp{W_0(\tfrac{1}{\delta}\me^{1/\delta})}\quad\text{for all $i$.} 
    \]
    Finally, using the property $W_0(x\me^x)=x$ for $x>0$, we conclude the proof. \qedhere
\end{enumerate}
\end{proof}

Lemma~\ref{lem:T} entails that the diagonal elements of $\mat{P}$ bound the diagonal elements of the $J_1$-generators of $\mat{P}$.
\begin{cor}\label{cor:qii_bounds}
    Let $\mat{P}=(p_{ij})$ be a $n\times n$ stochastic matrix. Let $\Delta=\max\{p_{11},\dots,p_{nn}\}$ and $\delta=\min\{p_{11},\dots,p_{nn}\}$. Suppose $\delta>0$. Then, if $\mat{Q}=(q_{ij})$ is a $J_1$-generator of $\mat{P}$, we have 
    \[
    1-\frac{1}{\delta}\leq q_{ii}\leq 1-\frac{1}{\Delta}\quad\text{for all $i$.}
    \]
\end{cor}
\begin{proof}
    If $\mat{Q}=(q_{ij})$ is a $J_1$-generator of $\mat{P}$, we have by proposition~\ref{prop:cembed:Tfixpt} that the vector $(\theta_1,\dots,\theta_n)$, where $\theta_i=\me^{1-q_{ii}}$ for all $i$, is a fixed oint of $\mat{T}$. Applying lemma~\ref{lem:T}\ref{lem:T:fp_in_X}, we have $\me^{1/\Delta}\leq\theta_i\leq\me^{1/\delta}$ for all $i$, from which the result follows readily.
\end{proof}

By combining propositions \eqref{prop:q_ij}, \eqref{prop:cembed:Tfixpt} and \eqref{prop:fixpT_to_J1gen}, it turns out that there is a one-to-one correspondence between the possible $J_1$-generators of $\mat{P}$ and the fixed points of the vector function $\mat{T}$. Regarding these fixed points, we now prove the following important result.

\begin{thm}\label{thm:TuniqueFP}
    Let the stochastic $n\times n$ matrix $\mat{P}=(p_{ij})$ be such that $p_{ii}>0$ for all $i$. Then, the vector function $\mat{T}:\Rset_+^n\to \Rset_+^n$, defined as in \eqref{eq:T}, has a unique fixed point.
\end{thm}
\begin{proof}
     From lemma~\ref{lem:T}\ref{lem:T:XtoX}, we know that $\mat{T}$ maps the compact convex set $\mathcal{X}\subset\Rset_+^n$, defined by \eqref{eq:X}, into itself. Also, $\mat{T}$ is continuous as the function $\rho$, defined by \eqref{eq:rho(x,y)}, is continuous (lemma~\ref{lem:rho}\ref{lem:rho:cont}) and continuity is preserved by linear combination and composition of continuous functions. Hence, by the Brouwer fixed-point theorem, $\mat{T}$ has a fixed point. By definition of $\mat{T}$, this fixed point must have all positive components. We now show that the function $\mat{g}:\Rset_+^n\to\Rset^n$ defined as $\mat{g}=\mat{T}-\mat{Id}$, where $\mat{Id}:\Rset_+^n\to\Rset_+^n$ is the identity mapping, satisfies all conditions of Theorem~3.1 in \cite{kennan2001uniqueness}. This theorem states sufficient conditions in order for the function $\mat{g}$ to have at most one vector $\mat{x}\in\Rset_+^n$ with $\mat{g}(\mat{x})=\mat{o}$. These conditions are (a) $\mat{g}$ is quasi-increasing and (b) $\mat{g}$ is strictly $R$-concave. Both (a) and (b) are proven in this paper, see lemma~\ref{lem:g.conditions}. So, we have established that $\mat{T}$ has exactly one fixed point.
\end{proof}

We are now in the position to formulate and prove our main theorem.

\begin{thm}\label{thm:main}
    Let the stochastic $n\times n$ matrix $\mat{P}=(p_{ij})$ be such that $p_{ii}>0$ for all $i$. Then, $\mat{P}$ has exactly one $J_1$-generator. Moreover, this $J_1$-generator $\mat{Q}=(q_{ij})$ has elements given by 
    \begin{equation}\label{eq:Q-fixp}
        q_{ii}= 1-\ln{\theta_i},\qquad q_{ij}=\frac{\rho(\theta_i,\theta_j)p_{ij}}{\theta_i p_{ii}}\quad(i\neq j)
    \end{equation}
    where the scalar function $\rho:\Rset_+^2\to\Rset_+$ is defined by \eqref{eq:rho(x,y)} and $(\theta_1,\dots,\theta_n)$ is the unique fixed point of the vector function $\mat{T}:\Rset_+^n\to\Rset_+^n$ defined  by \eqref{eq:T}.
\end{thm}
\begin{proof}
   We first prove that $\mat{P}$ has a $J_1$-generator. By theorem~\ref{thm:TuniqueFP}, the vector function $\mat{T}$ has a unique fixed point $\mat{\theta}=(\theta_1,\theta_2,\dots,\theta_n)\in\Rset_+^n$. Starting from $\mat{\theta}$, construct the matrix $\mat{Q}=(q_{ij})$ according to \eqref{eq:Q-fixp}. Then, $\mat{Q}$ is a $J_1$-generator of $\mat{P}$, by proposition~\ref{prop:fixpT_to_J1gen}. 
   
   To prove the uniqueness of the $J_1$-generator, suppose that $\mat{P}$ has  $J_1$-generators $\mat{R}=(r_{ij})$ and $\mat{S}=(s_{ij})$. Then, by proposition~\ref{prop:cembed:Tfixpt}, the vectors $\mat{\theta_{R}}=(\me^{1-r_{11}},\dots,\me^{1-r_{nn}})$ and $\mat{\theta_{S}}=(\me^{1-s_{11}},\dots,\me^{1-s_{nn}})$ are fixed points of the vector function $\mat{T}$. By theorem~\ref{thm:TuniqueFP}, $\mat{\theta_{R}}=\mat{\theta_{S}}$. Hence, by proposition~\ref{prop:q_ij}, we must have $\mat{R}=\mat{S}$. 
    
   Finally, the fact that a $J_1$-generator assumes the form \ref{eq:Q-fixp} is a consequence of proposition~\ref{prop:q_ij} and proposition~\ref{prop:cembed:Tfixpt}. The proof is now complete. 
   \end{proof}

According to theorem \ref{thm:main}, the unique $J_1$-generator is completely determined by the fixed point of the function $\mat{T}$. Besides, under the conditions of lemma \ref{lem:fixed point}, the function $\mat{T}$ is a contraction. Hence, under these conditions, an algorithm based on the fixed point iteration approach guarantees an appropriate estimation for the fixed point $(\theta_1,\dots,\theta_n)$ as outcome.

For a transition matrix $\mat{P}$ with identical positive diagonal elements, a closed-form formula for its unique $J_1$-generator can be given by virtue of corollary~\ref{cor:qii_bounds}. This type of transition matrices appear in e.g. models of DNA sequence evolution \cite{casanellas2020embeddability}.

\begin{cor} \label{closed form generator}
     Suppose $\mat{P}=(p_{ij})$ is a $n\times n$ stochastic matrix satisfying $p_{ii}=p>0$ for all $i$. Then, its unique $J_1$-generator $\mat{Q}$ is given by $\mat{Q}=\frac{1}{p}(\mat{P}-\mat{I})$, where $\mat{I}$ is the $n\times n$ identity matrix.
\end{cor}
\begin{proof}
    Let $\mat{Q}=(q_{ij})$ be a $J_1$-generator of $\mat{P}$. It follows from corollary~\ref{cor:qii_bounds} that $q_{ii}=1-1/p$ for all $i$. Moreover, if $i\neq j$, theorem~\ref{thm:main} and equation~\eqref{eq:rho(x,y)} imply that 
    \[
    q_{ij}=\frac{\rho(\me^{1-q_{ii}},\me^{1-q_{ii}})p_{ij}}{\me^{1-q_{ii}}p_{ii}}=\frac{\rho(\me^{1/p},\me^{1/p})p_{ij}}{\me^{1/p}p}=\frac{p_{ij}}{p}\text{.}
    \]
    In summary, we have
    \[
    q_{ii}=1-\frac{1}{p}=\frac{1}{p}(p_{ii}-1),\qquad
    q_{ij}=\frac{1}{p}p_{ij}\quad (i\neq j),
    \]
   concluding the proof.
\end{proof}

\section{Illustrations}

The aim of this section is twofold, namely (1) to illustrate the conditional embedding approach for some concrete transition matrices and (2) to compare the new approach with alternative low jump frequency approaches for embedding problems.
In \cite{jarrow1997markov} Jarrow, Lando and Turnbull found an approximation for the Markov generator in closed form under the model assumption that the probability of more than one jump per year is negligible. Their Markov generator $\mat{Q}_{\text{JLT}}=(q_{ij}^{\text{JLT}})$ is a product of the model assumption
\begin{equation}\label{eq:JLT:modelassumption}
    \Prob(H_i\geq 1\,|\,X_0=i)=p_{ii}\quad,\quad \Prob(H_i<1,X_{H_i}=j\,|\,X_0=i)= p_{ij}\quad(i\neq j)\text{,}
\end{equation}
where $H_i$ is the holding time in state $i$. This system of equations in the unknowns $q_{ij}^{\text{JLT}}$ can be solved explicitly to obtain
\begin{equation}\label{eq:JLT}
q^{\text{JLT}}_{ii}=\ln{p_{ii}} \quad,\quad q^{\text{JLT}}_{ij}=\frac{p_{ij} \ln{p_{ii}}}{p_{ii} - 1} \quad (i\neq j)\text{.}
\end{equation}

In this paper we study the embedding problem under the assumption 
\begin{equation}\label{empirical assumption}
    \Prob[X_1=j\,|\,X_0=i,N_{\!J}\leq 1]=p_{ij}\quad \text{for all $i$ and $j$.}
\end{equation}

In contrast to assumption~\eqref{eq:JLT:modelassumption}, assumption ~\eqref{empirical assumption} is about the data and not about the process. In fact, this paper does not preclude the model from having multiple transitions between $t=0$ and $t=1$. It does, however suppose that the data at hand are realisations of the underlying process with no more than one transition between $t=0$ and $t=1$. Hence, the estimated one step transition matrix is based solely on observations that did not jump more than once between $t=0$ and $t=1$.

    
Both approaches have the merit to make the identification phase as well as regularization redundant in the embedding problem. 



To avoid confusion, let us denote for the transition matrix $\mat{P}$, Jarrow's generator as $\mat{Q}_\text{JLT}$ and the $J_1$-generator as $\mat{Q}_{J_1}$. Then, an interesting question emerges of which of the matrices $\exp(\mat{Q}_\text{JLT})$ and $\exp(\mat{Q}_{J_1})$ is the best approximation to $\mat{P}$? The following section compares both approach for some interesting illustrations. As in \cite{davies2010embeddable}, throughout this paper the maximum absolute row sum is used as matrix norm, i.e. for an $n\times n$ matrix $\mat{M}=(m_{ij})$,  $||\mat{M}||_{\infty} = \max_{1 \leq i \leq n} \sum_{j=1}^{n} m_{ij}$.

\subsection{Credit rating transition matrix}
Consider the empirical transition matrix
\[
\mat{P}=
\begin{bmatrix}{}
  0.8910 & 0.0963 & 0.0078 & 0.0019 & 0.0030 & 0.0000 & 0.0000 & 0.0000 \\ 
  0.0086 & 0.9010 & 0.0747 & 0.0099 & 0.0029 & 0.0029 & 0.0000 & 0.0000 \\ 
  0.0009 & 0.0291 & 0.8896 & 0.0649 & 0.0101 & 0.0045 & 0.0000 & 0.0009 \\ 
  0.0006 & 0.0043 & 0.0656 & 0.8428 & 0.0644 & 0.0160 & 0.0018 & 0.0045 \\ 
  0.0004 & 0.0022 & 0.0079 & 0.0719 & 0.7765 & 0.1043 & 0.0127 & 0.0241 \\ 
  0.0000 & 0.0019 & 0.0031 & 0.0066 & 0.0517 & 0.8247 & 0.0435 & 0.0685 \\ 
  0.0000 & 0.0000 & 0.0116 & 0.0116 & 0.0203 & 0.0754 & 0.6492 & 0.2319 \\ 
  0.0000 & 0.0000 & 0.0000 & 0.0000 & 0.0000 & 0.0000 & 0.0000 & 1.0000 \\ 
  \end{bmatrix}
\]
based\footnote{We have adjusted five entries on the main diagonal to ensure all rows sum up to one.} on Table~3 in Jarrow et al. \cite[p.~506]{jarrow1997markov}. For this matrix, it can be proven that the vector function $\mat{T}:\Rset_+^{8}\to \Rset_+^{8}$, defined by \eqref{eq:T}, is a contraction mapping (according to \ref{lem:fixed point}). Using fixed-point iteration and \eqref{eq:Q-fixp}, we find that the unique $J_1$-generator, truncated to 4 decimal places, is:
\[
\mat{Q}_{J_1}=
\left[\begin{array}{@{}*{8}{S[table-format=-1.4]}@{\;}}
  -0.1221 & 0.1075 & 0.0088 & 0.0022 & 0.0036 & 0.0000 & 0.0000 & 0.0000 \\ 
  0.0096 & -0.1114 & 0.0836 & 0.0114 & 0.0035 & 0.0034 & 0.0000 & 0.0000 \\ 
  0.0010 & 0.0325 & -0.1271 & 0.0752 & 0.0122 & 0.0053 & 0.0000 & 0.0009 \\ 
  0.0007 & 0.0049 & 0.0755 & -0.1874 & 0.0798 & 0.0192 & 0.0024 & 0.0049 \\ 
  0.0005 & 0.0026 & 0.0094 & 0.0886 & -0.2759 & 0.1301 & 0.0178 & 0.0270 \\ 
  0.0000 & 0.0022 & 0.0036 & 0.0079 & 0.0647 & -0.2121 & 0.0592 & 0.0746 \\ 
  0.0000 & 0.0000 & 0.0152 & 0.0157 & 0.0287 & 0.1031 & -0.4460 & 0.2834 \\ 
  0.0000 & 0.0000 & 0.0000 & 0.0000 & 0.0000 & 0.0000 & 0.0000 & 0.0000 \\ 
  \end{array}
  \right]
\]

In contrast, the rate matrix $\mat{Q}_{\text{JLT}}$, published in Jarrow et al. \cite{jarrow1997markov} and defined by equation \eqref{eq:JLT}, is 
\[
\mat{Q}_{\text{JLT}}=\\
\left[\begin{array}{@{}*{8}{S[table-format=-1.4]}@{\;}}
  -0.1154 & 0.1020 & 0.0083 & 0.0020 & 0.0032 & 0.0000 & 0.0000 & 0.0000 \\ 
  0.0091 & -0.1043 & 0.0787 & 0.0104 & 0.0031 & 0.0031 & 0.0000 & 0.0000 \\ 
  0.0010 & 0.0308 & -0.1170 & 0.0688 & 0.0107 & 0.0048 & 0.0000 & 0.0010 \\ 
  0.0007 & 0.0047 & 0.0714 & -0.1710 & 0.0701 & 0.0174 & 0.0020 & 0.0049 \\ 
  0.0005 & 0.0025 & 0.0089 & 0.0814 & -0.2530 & 0.1180 & 0.0144 & 0.0273 \\ 
  0.0000 & 0.0021 & 0.0034 & 0.0073 & 0.0568 & -0.1927 & 0.0478 & 0.0753 \\ 
  0.0000 & 0.0000 & 0.0143 & 0.0143 & 0.0250 & 0.0929 & -0.4320 & 0.2856 \\ 
  0.0000 & 0.0000 & 0.0000 & 0.0000 & 0.0000 & 0.0000 & 0.0000 & 0.0000 \\ 
  \end{array}
  \right]
\]

Notice that the elements on the main diagonal of $\mat{Q}_{\text{JLT}}$ are in absolute value smaller than their counterparts in the $J_1$-generator. This property remains true for all matrices $\mat{P}$ having non-zero diagonal elements and follows from lemma~\ref{lem:J1_vs_QJLT_qii}.
It can be shown that $||\mat{P}-\exp{\mat{Q}_{J_1}}||_{\infty} < ||\mat{P}-\exp{\mat{Q}_{\text{JLT}}}||_{\infty}$.

\subsection{Transition matrices with same diagonal entries}
In case the stochastic matrix $\mat{P}$ has coinciding diagonal elements equal to $p>0$, according to corollary \ref{closed form generator} a closed-form solution to the equation $\mat{P}^{\{N_{\!J}\leq 1\}}=\mat{P}$ exists: $\mat{Q}_{J_1}=\frac{1}{p}(\mat{P}-\mat{I})$, where $\mat{I}$ is the $n\times n$ identity matrix. Besides, Jarrow's solution (see \cite[eqs. 30a \& 30b]{jarrow1997markov}) is $\mat{Q_{\text{JLT}}}=\frac{\ln{p}}{p-1}(\mat{P}-\mat{I})$. It is interesting to investigate to what extent $\exp\mat{Q}_{J_1}$ and $\exp\mat{Q_{\text{JLT}}}$ differ from $\mat{P}$. Note that both matrices $\mat{Q}_{J_1}$ and $\mat{Q_{\text{JLT}}}$ are of the form $\mat{Q}(k):=k(\mat{P}-\mat{I})$ with $k$ constant, and equal to $\frac{1}{p}$ and $\frac{\ln{p}}{p-1}$ respectively.

For the $2 \times 2$ case  $\mat{P}=\begin{bmatrix}
        p & 1-p \\
        1-p & p 
    \end{bmatrix}$ it is known that $\mat{P}$ is embeddable, with unique generator $\mat{Q}=\frac{\ln{(2p-1)}}{2(p-1)}(\mat{P}-\mat{I})$, if and only if $p > \frac{1}{2}$  (\cite{singer1974spilerman}). For transition matrices $\mat{P}$ with $p \leq \frac{1}{2}$ the conditional embedding approach results in a unique $J_1$-generator $\mat{Q}_{J_1}$ for which $\exp\mat{Q}_{J_1}$ is a better approximation to $\mat{P}$ than $\exp\mat{Q}_\text{JLT}$:

\begin{lem}
    Let 
    $\mat{P}=\begin{bmatrix}
        p & 1-p \\
        1-p & p 
    \end{bmatrix}$, where $0<p<1$ and define $\mat{Q}(k)=k(\mat{P}-\mat{I})$, where $\mat{I}$ is the $2\times 2$ identity matrix. Then, 
    \[
    ||\mat{P}-\exp{\mat{Q}(\tfrac{1}{p})}||_{\infty} < ||\mat{P}-\exp{\mat{Q}(\tfrac{\ln{p}}{p-1})}||_{\infty}\text{.}
    \]
\end{lem}
\begin{proof}
    It can be shown that 
    \[
    \mat{P}-\exp{\mat{Q}(k)}=
    \left(\frac{1}{2}+\frac{{\me}^{2 k \left(p-1\right)}}{2}-p\right)
\begin{bmatrix}
    -1 & 1 \\1 & -1
\end{bmatrix}\text{,}
    \]
    so that we have 
    \[
    ||\mat{P}-\exp{\mat{Q}(k)}||_{\infty}=|1+\me^{2k(p-1)}-2p|\text{.}
    \]
    Let $f(k)=1+\me^{2k(p-1)}-2p$. Since $p<1$, the function $f$ is strictly decreasing. Also, since $\ln{x}<x-1$ for all $x>0$ and $x\neq 1$, we have that $\ln{\frac{1}{p}}<\frac{1}{p}-1=\frac{1-p}{p}$ yielding $\frac{\ln{p}}{p-1}<\frac{1}{p}$. Hence, $f(\frac{\ln{p}}{p-1})>f(\frac{1}{p})$. Furthermore, $f(\frac{1}{p})>0$ ( lemma~\ref{lem:XXX}, 1.). Consequently,
    \[
    ||\mat{P}-\exp{\mat{Q}(\tfrac{1}{p})}||_{\infty}=|f(\tfrac{1}{p})|=f(\tfrac{1}{p})<f(\tfrac{\ln{p}}{p-1})=|f(\tfrac{\ln{p}}{p-1})|=||\mat{P}-\exp{\mat{Q}(\tfrac{\ln{p}}{p-1})}||_{\infty}\text{.}\qedhere
    \]
\end{proof}

Hence, for all $(2 \times 2)$ transition matrices with same diagonal entries, it is proven that $||\mat{P}-\exp{\mat{Q}_{J_1}}||_{\infty}<||\mat{P}-\exp{\mat{Q}_{\text{JLT}}}||_{\infty}$. For the $(3 \times 3)$ case, we investigate the transition matrices as introduced in lemma \ref{same diagonal entries}. Those transition matrices are not embeddable, since $p_{13}=0$ but $p^{(2)}_{13}=(1-p)^2/2>0$ \cite[Theorem~5, p.~126]{chung1967markov}. Since no generator does exist, it is worth to investigate the $J_1$-generator and Jarrow's generator. Lemma \ref{same diagonal entries} proves that $\exp\mat{Q}_{J_1}$ better approximates $\mat{P}$ than $\exp\mat{Q}_\text{JLT}$, i.e. $||\mat{P}-\exp{\mat{Q}_{J_1}}||_{\infty}<||\mat{P}-\exp{\mat{Q}_{\text{JLT}}}||_{\infty}$.

\begin{lem} \label{same diagonal entries}
Let $\mat{P}=
\begin{bmatrix}{}
  p & 1-p & 0 \\ 
  \frac{1}{2}(1-p) & p & \frac{1}{2}(1-p) \\ 
  0 & 1-p & p \\ 
  \end{bmatrix}
  $, where $0<p<1$ and define $\mat{Q}(k)=k(\mat{P}-\mat{I})$, where $\mat{I}$ is the $3\times 3$ identity matrix. Then, 
    \[
    ||\mat{P}-\exp{\mat{Q}(\tfrac{1}{p})}||_{\infty}<||\mat{P}-\exp{\mat{Q}(\tfrac{\ln{p}}{p-1})}||_{\infty}\text{.}
    \]
\end{lem}
\begin{proof}
It can be shown (e.g. using Sylvester's theorem for computing  functions of a matrix) that 
\[
\mat{P}-\exp{\mat{Q}}(k)=
\begin{bmatrix}
    -\alpha(k) & \beta(k) & \alpha(k)-\beta(k) \\[3pt] \frac{1}{2}\beta(k) & -\beta & \frac{1}{2}\beta(k) \\[3pt] \alpha(k)-\beta(k) & \beta(k) & -\alpha(k)
\end{bmatrix}
\]
where
\begin{equation*}\label{eq:ex:alphabeta}
    \alpha(k)=\tfrac{1}{4}\,\me^{2k(p-1)}+\tfrac{1}{2}\,\me^{k(p-1)}+\tfrac{1}{4}-p \quad\text{and}\quad \beta(k)=\tfrac{1}{2}\,\me^{2k(p-1)}+\tfrac{1}{2}-p\text{.}
\end{equation*}
It holds that
\begin{equation}\label{eq:alfabeta:1}
    \alpha(k)-\beta(k)=-\tfrac{1}{4}(1-\me^{k(p-1)})^2\leq 0\text{,}
\end{equation}
whence,
\begin{align}
    ||\mat{P}-\exp{\mat{Q}(k)}||_{\infty}
    &=|\alpha(k)|+|\beta(k)|+|\alpha(k)-\beta(k)| \nonumber\\
    &=|\alpha(k)|+|\beta(k)|+\beta(k)-\alpha(k)\text{.} \label{eq:norm:P-expQ(k)}
\end{align}
Note that $\alpha$ and $\beta$ are both strictly decreasing in $k$ since $p<1$. Also, we have 
\begin{align}
    \alpha(\tfrac{\ln{p}}{p-1})&=\tfrac{1}{4}(1-p)^2>0\text{,}\label{eq:alfabeta:2}\\
    \beta(\tfrac{\ln{p}}{p-1})&=\tfrac{1}{2}(1-p)^2>0\text{,}\label{eq:alfabeta:3}\\
    \beta(\tfrac{1}{p})&=\tfrac{1}{2}(\me^{2-2/p}+1-2p)>0\quad\text{by lemma~\ref{lem:XXX} (1.)}. \label{eq:alfabeta:4}
\end{align}
Hence, by \eqref{eq:norm:P-expQ(k)}, \eqref{eq:alfabeta:2} and \eqref{eq:alfabeta:3},
\begin{equation}\label{eq:norm:P-expQJ:3x3}
||\mat{P}-\exp{\mat{Q}(\tfrac{\ln{p}}{p-1})}||_{\infty}=
\alpha(\tfrac{\ln{p}}{p-1})+\beta(\tfrac{\ln{p}}{p-1})+\beta(\tfrac{\ln{p}}{p-1})-\alpha(\tfrac{\ln{p}}{p-1})=2\beta(\tfrac{\ln{p}}{p-1})\text{.}
\end{equation}

In case $\alpha(\frac{1}{p})\geq 0$, then by \eqref{eq:norm:P-expQ(k)} and \eqref{eq:alfabeta:4},
\[
||\mat{P}-\exp{\mat{Q}(\tfrac{1}{p})}||_{\infty}=\alpha(\tfrac{1}{p})+\beta(\tfrac{1}{p})+\beta(\tfrac{1}{p})-\alpha(\tfrac{1}{p})=2\beta(\tfrac{1}{p})\text{,}
\]
which yields $||\mat{P}-\exp{\mat{Q}(\tfrac{1}{p})}||_{\infty}<||\mat{P}-\exp{\mat{Q}(\tfrac{\ln{p}}{p-1})}||_{\infty}$, using \eqref{eq:norm:P-expQJ:3x3} and because $\frac{\ln{p}}{p-1}<\frac{1}{p}$ and $\beta$ is strictly decreasing.

In case $\alpha(\frac{1}{p})< 0$, we have by \eqref{eq:norm:P-expQ(k)} and \eqref{eq:alfabeta:1},
\begin{equation*}
||\mat{P}-\exp\mat{Q}(\tfrac{1}{p})||_{\infty}=-\alpha(\tfrac{1}{p})+\beta(\tfrac{1}{p})+\beta(\tfrac{1}{p})-\alpha(\tfrac{1}{p})= \tfrac{1}{2}(1-\me^{1-1/p})^2\text{.}
\end{equation*}
By lemma~\ref{lem:XXX} (2.) and the assumption $0<p<1$, it holds that $0<1-\me^{1-1/p}<\frac{4}{3}(1-p)$. Hence, using \eqref{eq:alfabeta:3} and \eqref{eq:norm:P-expQJ:3x3},
\[
||\mat{P}-\exp\mat{Q}(\tfrac{1}{p})||_{\infty}<\tfrac{8}{9}(1-p)^2<(1-p)^2=2\beta(\tfrac{\ln{p}}{p-1})=||\mat{P}-\exp\mat{Q}(\tfrac{\ln{p}}{p-1})||_{\infty}\text{.}
\]
In either case, we have proven the result.    
\end{proof}

\section*{Appendix: Lemma's and proofs}

\begin{lem}\label{lem:expWO.concave}
    The function $f$ with $f(t)=\me^{W_0(t)}$, $t\geq0$, is strictly concave.
\end{lem}
\begin{proof}
    By taking second order derivatives and since $W_0'(t)=\frac{W_0(t)}{t(1+W_0(t))}$ and $W_0''(t)=\frac{-2W_0(t)^2-W_0(t)^3}{t^2(1+W_0(t))^3}$ (see e.g. \cite{dence2013brief}), we find
    \[
f''(t)=f(t)\bigl(W_0'(t)^2+W_0''(t)\bigr)=f(t)\frac{-W_0(t)^2}{t^2(1+W_0(t))^3}\text{,}
    \]
    which is negative for all $t>0$ since $W_0(t)>0$ if $t>0$.
\end{proof}

Let $\mat{o}=(0,\dots,0)\in\Rset^n$. In what follows, we consider the partial ordering of $\Rset^n$ induced by componentwise ordering. For example, if $\mat{x}=(x_1,\dots,x_n)\in\Rset^n,$ and $\mat{y}=(y_1,\dots,y_n)\in\Rset^n$, we write $\mat{x}\preceq\mat{y}$ if and only if $x_i\leq y_i$ for all $i$. Likewise, we write $\mat{x}\succ\mat{o}$ if and only if $x_i>0$ for all $i$.

\begin{lem}\label{lem:rho}
    The function $\rho:\Rset_+^2\to\Rset_+$, defined as in \eqref{eq:rho(x,y)} , satisfies the following properties:
    \begin{enumerate}[label=(\arabic*)]
        \item It is linearly homogeneous, i.e. $\rho(\lambda\mat{x})=\lambda\rho(\mat{x})$ for all $\mat{x}\in\Rset_+^2$ and $\lambda>0$. \label{lem:rho:homog}
        \item It is continuous on $\Rset_+^2$. \label{lem:rho:cont}
        \item It is increasing, i.e. $\rho(\mat{x})\leq\rho(\mat{y})$ for all $\mat{x},\mat{y}\in\Rset_+^2$ with $\mat{x}\preceq\mat{y}$. \label{lem:rho:incr}
        \item $\min\{x,y\}\leq \rho(x,y)\leq\max\{x,y\}$ for all $(x,y)\in\Rset_+^2$. \label{lem:rho:maxmin}
   \end{enumerate}
\end{lem}
\begin{proof}
    Let $\mat{u}=(u_1,u_2)\in\Rset_+^2$. It is easy to see that $\rho(\mat{u})=u_2f(u_1/u_2)=u_1f(u_2/u_1)$, where $f$ is the continuous function defined by
    \[
    f(t)=\begin{cases}
        \frac{t\ln{t}}{t-1} & \text{if $t>0$ and $t\neq 1$}\\
        1 & \text{if $t=1$.}
    \end{cases}
    \]
    \begin{enumerate}[label=(\arabic*)]
        \item 
        Follows directly from the above.
        
        \item 
        A direct consequence of the above.
        
        \item 
        By standard calculus, we have 
        \[
        f'(t)=\begin{cases}
            \frac{t-1-\ln{t}}{(t-1)^2} & \text{if $t>0$ and $t\neq 1$} \\
            1/2 & \text{if $t=1$,}
        \end{cases}
        \]
        hence $f$ is (strictly) increasing on the positive real half-line since $\ln{t}<t-1$ for all $t>0$ with $t\neq 1$. Now, take $\mat{x}=(x_1,x_2)$ and $\mat{y}=(y_1,y_2)$, so that $\mat{o}\prec\mat{x}\preceq\mat{y}$. Then, as $f$ is increasing,
        $\rho(\mat{x})=x_2f(x_1/x_2)\leq x_2f(y_1/x_2)=y_1f(x_2/y_1)\leq y_1f(y_2/y_1)=\rho(\mat{y})$.

        \item 
        Consider the case $0<x\leq y$. Since $f$ is increasing on $\Rset_+$, we have 
        \[
        x=xf(1)\leq xf(y/x)=\rho(x,y)=yf(x/y)\leq yf(1)=y\text{.}
        \]
        The result for the case $0<y\leq x$ is proven analogously. \qedhere
    \end{enumerate}
\end{proof}

\begin{lem}\label{lem:T.increasing}
    The vector function $\mat{T}:\Rset_+^n\to\Rset_+^n$ from proposition~\ref{prop:cembed:Tfixpt} is increasing, i.e. $\mat{T}(\mat{x})\preceq\mat{T}(\mat{y})$ for all $\mat{x},\mat{y}\in\Rset_+^n$ with $\mat{x}\preceq\mat{y}$.
\end{lem}
\begin{proof}
    Let $i\in\{1,\dots,n\}$ and take $\mat{x},\mat{y}\in\Rset_+^n$ so that $\mat{x}\preceq\mat{y}$. Denote, for $\mat{u}=(u_1,\dots,u_n)$, $F_i(\mat{u})=\frac{1}{p_{ii}}\left(\sum_{j=1}^n p_{ij}\rho(u_i,u_j)\right)$. By lemma~\ref{lem:rho}\ref{lem:rho:incr}, we have $F_i(\mat{x})\leq F_i(\mat{y})$, which yields $T_i(\mat{x})=\me^{W_0(F_i(\mat{x}))}\leq \me^{W_0(F_i(\mat{y}))}=T_i(\mat{y})$ since the principal branch $W_0$ of the Lambert W function is increasing (see e.g. \cite{corless1996lambertw}).
\end{proof}

\begin{lem}\label{lem:g.conditions}
    Let the vector function $\mat{g}:\Rset_+^n\to\Rset^n$ given by $\mat{g}(\mat{x})=\mat{T}(\mat{x})-\mat{x}$ for all $\mat{x}\succ\mat{o}$, where the vector function $\mat{T}:\Rset_+^n\to\Rset_+^n$ is defined as in \eqref{eq:T}. Then, $\mat{g}=(g_1,\dots,g_n)$ is
    \begin{enumerate}[label=(\arabic*)]
        \item \emph{quasi-increasing}, i.e. for all $i$ holds that $\mat{o}\prec\mat{x}\preceq\mat{y}$ and $x_i=y_i$ imply $g_i(\mat{x})\leq g_i(\mat{y})$,
        \item \emph{strictly $R$-concave}, i.e. if $\mat{x}\succ\mat{o}$ and $\mat{g}(\mat{x})=\mat{o}$ and $0<\lambda<1$, then $\mat{g}(\lambda\mat{x})\succ\mat{o}$.
    \end{enumerate} 
\end{lem}
\begin{proof}$ $
\begin{enumerate}[label=(\arabic*)]
    \item 
    Take $i\in\{1,\dots,n\}$ and suppose $\mat{o}\prec\mat{x}\preceq\mat{y}$ with $x_i=y_i$. Then, $g_i(\mat{x})\leq g_i(\mat{y})$, since $x_i=y_i$ and $T_i(\mat{x})\leq T_i(\mat{y})$ (lemma~\ref{lem:T.increasing}). 

    \item 
    Let $\mat{x}=(x_1,\dots,x_n)\succ\mat{o}$ so that $\mat{g}(\mat{x})=\mat{o}$. Let $0<\lambda<1$ and take $i\in\{1,\dots,n\}$. Denote $F_i(\mat{x})=\frac{1}{p_{ii}}\left(\sum_{j=1}^n p_{ij}\rho(x_i,x_j)\right)$. By lemma~\ref{lem:rho}\ref{lem:rho:homog}, $F_i(\lambda\mat{x})=\lambda F_i(\mat{x})$. Hence,
    \[
    T_i(\lambda\mat{x})=\me^{W_0(F_i(\lambda\mat{x}))}=\me^{W_0(\lambda F_i(\mat{x}))}>\lambda\me^{W_0(F_i(\mat{x}))}=\lambda T_i(\mat{x})\text{,}
    \]
    where the inequality follows from the fact that the function $t\mapsto\me^{W_0(t)}$ is strictly concave (lemma~\ref{lem:expWO.concave}) and $W_0(0)=0$. Consequently, for all $i$,
    \[
    g_i(\lambda\mat{x})=T_i(\lambda\mat{x})-\lambda x_i>\lambda(T_i(\mat{x})-x_i)=\lambda g_i(\mat{x})=0\text{,}
    \]
    i.e. $\mat{g}(\lambda\mat{x})\succ\mat{o}$. \qedhere
\end{enumerate}
\end{proof}

\begin{lem}\label{lem:J1_vs_QJLT_qii}
    Let the stochastic $n\times n$ matrix $\mat{P}=(p_{ij})$ be such that $p_{ii}>0$ for all $i$. Let $\mat{Q}=(q_{ij})$ be the $J_1$-generator of $\mat{P}$. Then, $q_{ii}\leq\ln{p_{ii}}$ for all $i$.
\end{lem}
\begin{proof}
    By theorem~\ref{thm:main}, $q_{ii}=1-\ln{\theta_i}$ for all $i$, where $\theta=(\theta_1,\dots,\theta_n)$ is the unique fixed point of the vector function $\mat{T}=(T_1,\dots,T_n)$, defined by \eqref{eq:T}. Take $i\in\{1,\dots,n\}$. Thus, 
    \[
    \theta_i=T_i(\theta)=\exp{W_0(\tfrac{1}{p_{ii}}\sum_{j}p_{ij}\rho(\theta_i,\theta_j))}\text{,}
    \]
    which yields
    \[
    \theta_i\ln{\theta_i}=\frac{1}{p_{ii}}\sum_{j}p_{ij}\rho(\theta_i,\theta_j)\text{.}
    \]
    Using the definition of the function $\rho$ in \eqref{eq:rho(x,y)}, the above equation can be rewritten as 
    \begin{equation}\label{eq:theta_i_e}
    \left(\frac{\theta_i}{\me}-1\right)\rho(\theta_i,\me)=\theta_i\ln{\theta_i}-\theta_i= \frac{1}{p_{ii}}\sum_{j:j\neq i}p_{ij}\rho(\theta_i,\theta_j)\text{.}
    \end{equation}
    Now, since $1-\ln{\theta_j}=q_{jj}\leq 0$, we have that $\theta_j\geq\me$ for all $j$. By lemma~\ref{lem:rho}\ref{lem:rho:incr},  $\rho(\theta_i,\theta_j)\geq \rho(\theta_i,\me)$ for all $j$. Hence, it follows from \eqref{eq:theta_i_e} that
    \[
    \left(\frac{\theta_i}{\me}-1\right)\rho(\theta_i,\me)\geq
    \frac{1}{p_{ii}}\sum_{j:j\neq i}p_{ij}\rho(\theta_i,\me)=\frac{1-p_{ii}}{p_{ii}}\rho(\theta_i,\me)\text{,}
    \]
    which simplifies to $\theta_i\geq \me/p_{ii}$. Upon taking logarithms of both sides of this inequality and using $q_{ii}=1-\ln{\theta_i}$, we arrive at $q_{ii}\leq\ln{p_{ii}}$. 
\end{proof}

\begin{lem}\label{lem:fixed point}
Let $\mat{P}=(p_{ij})$ be a $n\times n$ stochastic matrix. Let $\Delta=\max\{p_{11},\dots,p_{nn}\}$ and $\delta=\min\{p_{11},\dots,p_{nn}\}$. Suppose $\delta>0$. Then, with respect to the max norm $||\cdot||_{\infty}$, the function $\mat{T}:\mathcal{X}\to \mathcal{X}$, defined as in  \eqref{eq:T}, is Lipschitz continuous with Lipschitz constant $K=\frac{1+(\frac{1}{\delta}-1)C(\alpha)}{1+\frac{1}{\Delta}}$ where $C(\alpha)=-1+\frac{\alpha+1}{\alpha-1}\ln{\alpha}$ and $\alpha=\me^{\frac{1}{\delta}-\frac{1}{\Delta}}$.
\end{lem}

\begin{lem}\label{lem:XXX}
    For all $0<p<1$, it holds that 
    \begin{enumerate}
        \item $1+\me^{2-2/p}-2p>0$,
        \item $1-\me^{1-1/p}<\frac{4}{3}(1-p)$.
    \end{enumerate}
\end{lem}
\begin{proof}
    To prove the first inequality, let $f(t)=1+\me^{2-2/t}-2t$, which is continuous on the half-open interval $(0,1]$.  A straightforward calculation reveals that $f'(t)=2t^{-2}(\me^{2-2/t}-t^2)$ and $f''(t)=4t^{-4}(1-t)\me^{2-2/t}$. So, $f''(t)>0$ for all $t\in(0,1)$. Consequently, $f'(t)<0$ for all $t\in(0,1)$ because $f'$ is monotone increasing on $(0,1)$ and $f'(1)=0$. Hence, $f$ is monotone decreasing on $(0,1)$. The result now follows from the fact that $f(1)=0$.

    To prove the second inequality, consider the function $f(t)=\frac{4}{3}(1-t)-1+\me^{1-1/t}$ which is differentiable on $\{t\in\Rset\,|\,t>0\}$. Let $p_0$ be a critical point of $f$, then $f'(p_0)=-\frac{4}{3}+\me^{1-1/p_0}p_0^{-2}=0$, yielding $\me^{1-1/p_0}=\frac{4}{3}p_0^2$. Clearly, $p_0\neq\frac{1}{2}$, hence, 
    \[
    f(p_0)=\tfrac{4}{3}(1-p_0)-1+\tfrac{4}{3}p_0^2=\tfrac{4}{3}(\tfrac{1}{4}-p_0+p_0^2)=\tfrac{4}{3}(\tfrac{1}{2}-p_0)^2>0\text{.}
    \]
    So, all critical points of $f$ have positive function values. In addition, we have $f(1)=0$ and $\lim_{t\to0^+}f(t)=1/3>0$. Therefore, $f(t)>0$ for all $t\in(0,1)$. 
\end{proof}

\bibliographystyle{unsrt}
\bibliography{condem-ref}

\end{document}